\input amstex
\input amsppt.sty
\nologo

\def\Rls{\Bbb{R}}

\topmatter
\title Fully nonlinear parabolic equations in two space variables
\endtitle
\author Ben Andrews\endauthor
\affil Centre for Mathematics and its Applications\\
Australian National University\endaffil
\address Centre for Mathematics and its Applications, 
Australian National University, A.C.T. 0200, Australia
\endaddress
\email andrews\@maths.anu.edu.au
\endemail
\subjclass 35K55, 35B65
\endsubjclass
\abstract{H\"older estimates for spatial second derivatives are proved
for solutions of fully nonlinear parabolic equations in two space
variables.  Related techniques extend the regularity theory 
for fully nonlinear parabolic equations in higher dimensions.
}
\endabstract
\endtopmatter
\document
\rightheadtext{Fully nonlinear parabolic equations}

\head 1.  Introduction\endhead

Elliptic equations in two variables are very well understood, and
the regularity theory for such equations is significantly stronger
than that available for elliptic equations in higher dimensions.
In particular, Morrey \cite{M} and Nirenberg \cite{N} proved 
H\"older estimates for the first derivatives of
solutions of uniformly elliptic equations in two variables,
depending only on bounds for the coefficients:

\proclaim{Theorem 1}
Let $\Omega\subset\Rls^2$, and set
$d_{\Omega}(z)=d(z,\partial\Omega)$ for all $z\in\Omega$.  Let $u$ be a
bounded
$C^2(\Omega)$ solution of
$$
a u_{xx} + 2b u_{xy} + c u_{yy} = f\tag{1}
$$
where $a$, $b$ and $c$ are measureable functions on $\Omega$ with
$\bmatrix a &b\cr b &c\cr 
\endbmatrix\geq \lambda \sqrt{ac-b^2}I$ and $\delta=\sqrt{ac-b^2}>0$
everywhere.  Then for any $\alpha\in(0,\sqrt{\lambda})$ there exists
$C=C(\lambda,\alpha)$ such that for all points $p\neq q$ in $\Omega$
with $d=\min\{d_\Omega(p),d_\Omega(q)\}>0$,
$$
{|Du(q)-Du(p)|\over|p-q|^{\alpha}}\leq Cd^{-\alpha}
\sup_{\Omega}\left(|Du|+
d_{\Omega}\delta^{-1}|f|\right).
$$
\endproclaim

These estimates can also be applied to fully nonlinear uniformly 
elliptic equations in two variables, to give H\"older estimates for
second derivatives.  in this case the equations have the form
$$
F[u]=F(D^2u,Du,u,x)=0\tag{2}
$$
where $u:\Omega\subset\Rls^2\to\Rls$, and $F:
{\Bbb S}_2\times\Rls^2\times\Rls\times\Omega\to\Rls$ is Lipschitz
in all variables (here ${\Bbb S}_2$ is the space of symmetric
$2\times 2$ matrices) and
uniformly monotone in the first argument, so that there exist constants
$\Lambda\geq\lambda>0$ such that
$$
\lambda I\leq [\dot F^{ij}]\leq \Lambda I\tag{3}
$$
where $\dot F^{ij} = {\partial F(r,p,z,x)\over\partial r_{ij}}$.

\proclaim{Theorem 2}
Let $u\in C^3(\Omega)$ satisfy $F[u]=0$ in $\Omega\subset\Rls^2$,
and suppose \thetag{2} holds.  Then there exists
$\alpha=\alpha(\lambda/\Lambda)$ such that for any points
$p\neq q$ in
$\Omega$ with $d=\min\{d_\Omega(p),d_\Omega(q)\}>0$, 
$$
{|D^2u(p)-D^2u(q)|\over|p-q|^{\alpha}}\leq Cd^{-\alpha}
\sup_{\Omega}
\left(|D^2u|+|Du|+1\right)
$$
where $C$ depends only on $\lambda/\Lambda$ and $\sup|DF|$.
\endproclaim

In contrast, the situation in higher dimensions is much worse:  For
fully nonlinear equations there is no H\"older estimate known for
the second derivatives of solutions, unless the equation satisfies
a concavity condition with respect to the components of the
second derivatives.  The best result available is the following,
due to Evans \cite{E1--2} and Krylov \cite{Kr} (I follow the treatment
in \cite{GT}):

\proclaim{Theorem 3}
Let $u\in C^4(\Omega)$ satisfy $F[u]=0$ in $\Omega$ where $F$ is a $C^2$
function of the form \thetag{2} which is uniformly elliptic (so that
$\lambda I\leq\dot F^{ij}\leq\Lambda I$ for some $\Lambda\geq\lambda>0$)
and concave with respect to the first argument.  Then for any
$\Omega'\subset\subset\Omega$,
$$
\sup_{p,q\in\Omega'}{|D^2u(p)-D^2u(q)|\over |p-q|^\alpha}\leq C
$$
where $\alpha$ depends on $n$, $\lambda$ and $\Lambda$, and
$C$ depends on $n$, $\lambda$, $\Lambda$, $|u|_{C^2(\Omega)}$
and $d(\Omega',\partial\Omega)$, and on bounds for the first
and second derivatives of $F$ (other than the second derivative in the
first argument).
\endproclaim

In the parabolic case, there are results similar to Theorem 3 (due to
Krylov \cite{Kr}):

\proclaim{Theorem 4}
Let $u\in C^4(\Omega\times(0,T])$ satisfy
$$
{\partial u\over\partial t} = F(D^2u,Du,u,x,t)
$$
where $F$ is $C^2$, $\lambda I\leq [\dot F^{ij}]\leq \Lambda I$
for some $0<\lambda\leq\Lambda$, and $F$ is concave in the first
argument.  Then for any $\tau>0$ and $\Omega'\subset\subset\Omega$, 
$$
\eqalign{
\sup_{s,t\in[\tau,T],p,q\in\Omega'}&
\left({|D^2u(p,t)-D^2u(q,t)|\over|p-q|^{\alpha}+|s-t|^{\alpha/2}}
+{|\partial_tu(p,t)-\partial_tu(q,t)|\over|p-q|^{\alpha}+|s-t|^{\alpha/2}}\right)\cr
&\null+\sup_{p\in\Omega', \tau\leq s, t\leq T}
{|Du(p,t)-Du(p,s)|\over |s-t|^{(1+\alpha)/2}}\qquad\leq C\cr
}
$$
where $\alpha$ depends on $n$, $\lambda$ and $\Lambda$, and
$C$ depends on $n$, $\lambda$, $\Lambda$,
$\sup_{\Omega\times(0,T]}|D^2u|$, $\sup_{\Omega\times(0,T]}|\partial_tu|$,
$d(\Omega',\partial\Omega)$, $\tau$ and bounds for the first and second
derivatives of $F$ (other than the second derivative in the first
argument).
\endproclaim

The Morrey and 
Nirenberg estimates rely either on quasiconformal mapping estimates
or on the fact that in two dimensions the first derivatives of solutions
satisfy divergence-form elliptic equations.  Neither of these
methods seems to generalise readily to the parabolic setting.  There
are, however, special estimates known for parabolic equations in one space
variable, due largely to Kruzhkov \cite{Kz}.

In this paper I will prove an analogue of Theorem 2 for parabolic
equations in two space variables.  I do not know whether an
analogue of Theorem 1 holds.

\proclaim{Theorem 5} Let $\Omega$ be a domain in $\Rls^2$.
Let $u\in C^3(\Omega\times(0,T])$ be a solution of the fully nonlinear
equation 
$$
{\partial u\over\partial t}= F(D^2u,Du,u,x,t)
$$
where $F$ is Lipschitz in all arguments and uniformly monotone in the
first argument, so that $\lambda I\leq [\dot F^{ij}]\leq \Lambda I$
for some $\Lambda\geq \lambda>0$.  Then for any $\tau\in(0,T)$ and 
$\Omega'\subset\subset \Omega$, 
$$
\eqalign{
\sup_{s,t\in[\tau,T],p,q\in\Omega'}&
\left({|D^2u(p,t)-D^2u(q,t)|\over|p-q|^{\alpha}+|s-t|^{\alpha/2}}
+{|\partial_tu(p,t)-\partial_tu(q,t)|\over|p-q|^{\alpha}+|s-t|^{\alpha/2}}\right)\cr
&\null+\sup_{p\in\Omega', \tau\leq s, t\leq T}
{|Du(p,t)-Du(p,s)|\over |s-t|^{(1+\alpha)/2}}\qquad\leq C\cr
}
$$
where $\alpha$ depends on $\lambda$ and $\Lambda$, and
$C$ depends on $\lambda$, $\Lambda$,
$\sup_{\Omega\times(0,T]}(|D^2u|+|Du|)$,
$\sup_{\Omega\times(0,T]}|\partial_tu|$,
$d(\Omega',\partial\Omega)$, $\tau$, and bounds for the first
derivatives of $F$.
\endproclaim

This gives a result of similar strength to Theorem 2 for fully nonlinear
parabolic equations in two space variables.

I will also apply similar ideas to parabolic equations in higher
dimensions, to relax the requirement of concavity to allow just convexity
of level sets.  

\proclaim{Theorem 6} Let $\Omega$ be a domain in $\Rls^n$.
Suppose $u\in C^4(\Omega\times(0,T])$ satisfies
$$
{\partial u\over\partial t}=F(D^2u,Du,u,x,t)
$$
where $F$ is $C^2$ and $\lambda I\leq [\dot F^{ij}]\leq
\Lambda I$ for some $\Lambda\geq\lambda>0$, and 
$\ddot F^{ij,kl}M_{ij}M_{kl}\leq 0$ for all
matrices
$[M_{ij}]$ for which $\dot F^{ij}M_{ij}=0$.  Then for any
$\Omega'\subset\subset\Omega$ and
$\tau\in(0,T)$,
$$
\eqalign{
\sup_{s,t\in[\tau,T],p,q\in\Omega'}&
\left({|D^2u(p,t)-D^2u(q,t)|\over|p-q|^{\alpha}+|s-t|^{\alpha/2}}
+{|\partial_tu(p,t)-\partial_tu(q,t)|\over|p-q|^{\alpha}+|s-t|^{\alpha/2}}\right)\cr
&\null+\sup_{x\in\Omega', \tau\leq s, t\leq T}
{|Du(x,t)-Du(x,s)|\over |s-t|^{(1+\alpha)/2}}\qquad\leq C\cr
}
$$
where $\alpha$ depends on $n$, $\lambda$ and $\Lambda$, and $C$
depends on $\lambda$, $\Lambda$,
$\sup_{\Omega\times(0,T]}(|D^2u|+|Du|)$,
$\sup_{\Omega\times(0,T]}|\partial_tu|$,
$d(\Omega',\partial\Omega)$, $\tau$, $K$, and bounds for first and
second derivatives of $F$.
\endproclaim

In the elliptic case, this extension
is trivial because the equation $F=0$ is the same as the equation
$\psi(F)=0$ if $\psi$ is an increasing function with $\psi(0)=0$. However, in
the parabolic case there is a big difference between the two equations
$$
\frac{\partial u}{\partial t}=F
$$
and 
$$
\frac{\partial u}{\partial t}=\psi(F).
$$
If $F$ is concave in $D^2u$ and $\psi$ is increasing, then the
latter equation is covered by Theorem 6 but not in general by Theorem 4.

I would like to thank Neil Trudinger and Craig Evans for useful
suggestions.

\head 2.  Some background results\endhead

In this section I will recall some results and notation to be used in the
proofs of Theorems 5 and 6.

\subhead 2.1 Function spaces\endsubhead

Let $\Omega$ be a domain in $\Rls^n$.  
As usual I will denote by $C^{k,\alpha}$ the space of functions which
have all derivatives up to order $k$ H\"older-continuous of with exponent
$\alpha\in (0,1]$.  

In the parabolic setting, where we are considering solutions on a space-time
region $Q=\Omega\times(0,T]$, it is useful to introduce the spaces
$P^{k,\alpha}(Q)$ consisting of functions $u$ on $Q$ which have the
following norm bounded:
$$
\eqalign{
|u|_{P^{k,\alpha}(Q)} &= \sum_{|a|+2b\leq k}\sup_{(x,t)\in Q}
|\partial_t^bD^au(x,t)|\cr
&\quad\null+\sum_{|a|+2b=k}\sup_{(x,t)\neq (y,s)\in Q}
{|\partial_t^bD^au(x,t)-\partial_t^bD^au(y,s)|
\over |x-y|^\alpha+|s-t|^{\alpha\over 2}}\cr
&\quad\null+\sum_{|a|+2b=k-1}\sup_{(x,t)\neq(x,s)\in Q}
{|\partial_t^bD^au(x,t)-\partial_t^bD^au(x,s)|\over |s-t|^{{1+\alpha\over 2}}}
\cr}
$$

\subhead 2.2 The Krylov-Safonov H\"older estimate\endsubhead

The Krylov-Safonov Harnack inequality, first proved in \cite{KS}, provides 
an oscillation  estimate for solutions of elliptic or parabolic equations.  
The following H\"older estimate is a consequence of this.   

\proclaim{Theorem 7}
Denote by $Q_r$
the region $B_{r}\times[-r^2,0]\subset{\Bbb R}^n\times{\Bbb R}$.
Let $u: Q_R\to{\Bbb R}$ be a smooth solution of an equation
of the form
$$
\frac{\partial u}{\partial t} = a^{ij}D_iD_ju+b^iD_iu+cu+f
$$
where the coefficients are bounded and measureable, and $\lambda
I\leq [a^{ij}]\leq\Lambda I$.  Then for some $C>0$ and $\alpha\in(0,1)$
depending only on $n$,
$\lambda/\Lambda$, and bounds for coefficients,
$$
|u|_{P^{0,\alpha}(Q_{R/2})}
\leq C\left(|u|_{P^0(Q_R)}+|f|_{L^\infty(Q_R)}\right).
$$
\endproclaim

\subhead 2.3 A characterisation of $C^{1,\alpha}$ functions\endsubhead

I will make use of the following characterisation of $C^{1,\alpha}$
functions in terms of difference quotients (see \cite{T}):

\proclaim{Theorem 8}
Let $u$ be a smooth function on $B_R\subset{\Bbb R}^n$, and suppose there exists a constant $C_0$ such that for all $\xi\in S^{n-1}$, $h>0$, and $x\in B_{R-h}$,
$$
|u(x+h\xi)+u(x-h\xi)-2u(x)|\leq Ch^{1+\alpha},
$$
where $\alpha\in(0,1]$.  Then 
$$
\sup_{x\neq y\in B_{R/2}}\frac{|Du(x)-Du(y)|}{|x-y|^\alpha}\leq C(\alpha)C_0.
$$
\endproclaim

\head 3. Two space variables\endhead

Assume (by rescaling if necessary) that $u$ is a $C^3$ solution on the domain 
$Q_1 = B_1(0)\times (-1,0]\subset\Rls^2\times\Rls$ of
a fully nonlinear parabolic equation
$$
{\partial u\over\partial t}=F(D^2u,Du,u,x,t)
$$
where $F$ is defined on a convex set ${\Cal S}$ in ${\Bbb
S}_2\times\Rls^2\times\Rls
\times B_1\times(-1,0]$ containing
$\{{\Cal J}u=(D^2u(x,t), Du(x,t), u(x,t), x,t): (x,t)\in Q_1\}$.
Assume that $F$ is Lipschitz on ${\Cal S}$, and satisfies
$$
\lambda \text{\rm I}\leq [\dot F^{ij}]\leq\Lambda \text{\rm I}
$$
at each point of ${\Cal S}$, for some $\Lambda\geq\lambda>0$.  We will obtain
estimates on regions $Q_r = B_r(0)\times (-r^2,0]$ for suitably small $r$.

\subhead 3.1 Regularity of the time derivative\endsubhead

Let $\tau\in(0,1)$, and define
$v_\tau(x,t)={1\over\tau}\left(u(x,t)-u(x,t-\tau)\right)$.  Then $v_{\tau}$
satisfies a parabolic equation on $B_1\times (\tau-1,0]$:
$$
\eqalign{
{\partial v_\tau\over\partial t} &= {
F({\Cal J}u(x,t))
-F({\Cal J}u(x,t-\tau))
\over\tau}\cr
&= {1\over\tau}\int_0^1DF|_{s{\Cal J}u(x,t)+(1-s){\Cal J}u(x,t-\tau)}\cdot
\left({\Cal J}u(x,t)-{\Cal J}u(x,t-\tau)\right)\,ds\cr
&=a^{ij}D_iD_jv_{\tau}+b^iD_iv_{\tau}+cv_\tau+f.
\cr
}
$$
Here, writing ${\Cal J}(s)=s{\Cal J}u(x,t)+(1-s){\Cal J}u(x,t-\tau)$, the
coefficients are given by
$a^{ij} = \int_0^1\dot F^{ij}|_{{\Cal J}(s)}\,ds$, $b^i=\int_0^1
F^{p^i}|_{{\Cal J}(s)}\,ds$, $c=\int_0^1F^z|_{{\Cal J}(s)}\,ds$ and
$f=\int_0^1{\partial F\over\partial t}|_{{\Cal J}(s)}\,ds$.  Note that
$\lambda \text{\rm I}\leq a^{ij}\leq \Lambda\text{\rm I}$, and
$|b^i|+|c|+|f|\leq C\text{\rm Lip}(F)$.

Theorem 7 applies to give an oscillation
estimate for $v_\tau$ independent of $\tau$, and hence also for $u_t$.
In particular, $u_t$ is H\"older continuous in both space and time, and for
$(x_1,t_1)$ and $(x_2,t_2)$ in $Q_{1/2}$
$$
|u_t(x_1,t_1)-u_t(x_2,t_2)|\leq C\left(|x_2-x_1|^2+|t_2-t_1|
\right)^{\alpha/2}\left(\left|{\partial F\over\partial
t}\right|_{L^\infty(Q_1)} +\left|u_t\right|_{L^\infty(Q^1)}
\right).
$$

\subhead 3.2 Spatial regularity of second space derivatives\endsubhead

This is the key estimate:  When restricted to each time slice, the function
$u$ has H\"older-continuous second derivatives.  This follows by combining
the Morrey-Nirenberg estimates for elliptic equations with
either a perturbation argument or an argument using difference quotients.  I 
present the latter argument.

Fix $t\in (-1/4,0]$.  Write $w(x)=u(x,t)$,
and $\phi(x)=u_t(x,t)$.  Then the following elliptic equation holds:
$$
G[D^2w(x),Dw(x),w(x),x]=\phi(x),
$$
where $G[r,p,z,x]=F[r,p,z,x,t]$ and $\phi(x)=u_t(x,t)$.  Note that
$\lambda\text{\rm I}\leq
\dot G\leq
\Lambda\text{\rm I}$, 
$\text{\rm Lip}(G)\leq\text{\rm Lip}(F)$, and $\phi$ satisfies the oscillation
estimate derived in Section 3.1.

Let $\xi$ be a unit vector, and for $h>0$ let $\delta_hw$ be the
difference quotient of $w$ in the direction $\xi$ with step $h$:
$$
\delta_hw(x) = h^{-1}\left(w(x+h\xi)-w(x)\right).
$$
$\delta_hw$ is defined on $B_{1-h}(0)$, and
satisfies an elliptic equation:  Denote ${\Cal
D}w(x)=(D^2w(x),Dw(x),w(x),x)$.  Then
$$
\eqalign{
0&=h^{-1}\left(G[{\Cal D}w(x+h\xi)]-G[{\Cal D}w(x)]\right)-\delta_h\phi\cr
&=h^{-1}\int_0^1DG|_{{\Cal D}(s)}\cdot
\left({\Cal D}w(x+\xi h)-{\Cal D}w(x)\right)\,ds-\delta_h\phi\cr
&=\tilde a^{ij}D_iD_j\delta_hw+\varphi
\cr
}
$$
Here ${\Cal D}(s) = s{\Cal D}w(x+h\xi)+(1-s){\Cal D}w(x)$ and 
$\varphi=\tilde b^iD_i\delta_hw+\tilde c\delta_hw+\tilde f-\delta_h\phi$,
where $\tilde f = \xi^i\int_0^1 {\partial G\over\partial x^i}|_{{\Cal
D}(s)}\,ds$,
$\tilde a^{ij} = \int_0^1\dot G^{ij}|_{{\Cal D}(s)}\,ds$, 
$\tilde b^i = \int_0^1 G^{p^i}|_{{\Cal D}(s)}\,ds$, and
$\tilde c = \int_0^1 G^z|_{{\Cal D}(s)}\,ds$. 
It follows that $\lambda|v|^2\leq \tilde a^{ij}v_iv_j\leq \Lambda|v|^2$ for
all $v\neq 0$, and 
$|\varphi(x)|\leq C\text{\rm
Lip}(F)(1+|u|_{P^2(Q_r)})+Ch^{\alpha-1}(\text{\rm
Lip}(F)+|u|_{P^2(Q_1)})$ for $x\in B_{1/4}$ and $h<1/4$.

The Morrey-Nirenberg estimates of Theorem 1 apply to the above equation to
give the following estimate:  If $x\in B_{1/8}(0)$ and $h<r\leq 1/16$,
then
$$
|D\delta_hw(x+h\xi)-D\delta_hw(x)|\leq C(h/r)^{\alpha}
\left(|u|_{P^2(Q_R)}+r|\varphi|\right).
$$
where $C$ depends on $\Lambda/\lambda$.
In particular, if we choose $r=h^{\beta}$ for $\beta\in(0,1)$ then
we have for $x\in B_{1/8}$ and $h<2^{-4/\beta}$ 
$$
|\delta_h\delta_hDw|\leq Ch^{-1}(1+|u|_{P^2(Q_1)})
\left(h^{\alpha(1-\beta)}
+h^{\alpha-(1-\beta)(1-\alpha)}\right).
$$
If we choose $1>\beta>\max\{0,{1-2\alpha\over1-\alpha}\}$ then we find
$$
|\delta_h\delta_hDw|\leq
Ch^{-1+\varepsilon}(1+|u|_{P^2(Q_1)}),
$$
for some $\varepsilon\in(0,1]$ depending on $\alpha$ and $\beta$,
and it follows from the characterisation of $C^{1,\alpha}$ functions in
section 2.4 that the second derivatives of $w$ are H\"older-continuous
with exponent $\varepsilon$, and
$$
|w|_{C^{2,\varepsilon}(B_{1/8})}\leq C(1+|u|_{P^2(Q_1)})
$$
where $C$ depends on $\lambda$ and $\text{\rm Lip}(F)$.  This proves the
required spatial regularity of spatial second derivatives on the region 
$Q_{1/8}$.

\subhead 3.3 Time regularity of first space derivatives\endsubhead

The spatial $C^{2,\alpha}$ estimate established in the previous section
can be used to deduce an estimate on the continuity of first spatial
derivatives in time, using the parabolic maximum principle.

Let $\xi$ be a unit vector.  The function $\delta_hu$
defined by $\delta_hu(x)=(u(x+h\xi)-u(x))/h$ satisfies a useful
evolution equation:
$$
{\partial\over\partial t}\delta_hu = a^{ij}D_iD_j\delta_hu+\varphi
\tag{3.3.1}
$$
where $\sup_{Q_{1/2}}|\varphi|\leq C(\text{\rm Lip}(F), |u|_{P^2(Q_1)})$.
The result of Section 3.2 shows that $\delta_hu$ is $C^{1,\alpha}$ (uniformly
in $h$), so that
$$
\left|\delta_hu(z',t)-\delta_hu(z,t)-D\delta_hu(z,t)\cdot(z'-z)\right|\leq
C|z'-z|^{1+\alpha}
$$
on the region $Q_{1/8-h}$.  Young's inequality gives the estimate
$$
|z'-z|^{1+\alpha}\leq
\varepsilon+C\varepsilon^{-{1-\alpha\over 1+\alpha}}|z'-z|^2,
$$
for any $\varepsilon>0$, and therefore
$$
\delta_hu(z',t)\leq \delta_hu(z,t)+D\delta_hu(z,t)\cdot(z'-z)
+\varepsilon+C\varepsilon^{-{1-\alpha\over1+\alpha}}|z'-z|^2.
$$
By the bounds on $a^{ij}$ and $\varphi$, the function
$$
\eqalign{
\Psi_+(z',t')&=\delta_hu(z,t)+D\delta_hu(z,t)\cdot(z'-z)
+\varepsilon+C\varepsilon^{-{1-\alpha\over1+\alpha}}|z'-z|^2\cr
&\quad\null
+(\sup|\varphi|+4C\Lambda\varepsilon^{-{1-\alpha\over1+\alpha}})(t'-t)\cr}
$$
is a supersolution of Equation (3.3.1) on $Q_{1/16}$ provided $h<1/32$, and
if $z\in B_{1/32}$ and 
$\varepsilon<C(\Lambda,\alpha)|u|_{P^2(Q_1)}^{-{1+\alpha\over 1-\alpha}}$
then this supersolution is above
$\delta_hu$ on the boundary $\partial B_{1/16}\times (-1/256,0]$.  Therefore
by the parabolic maximum principle
$\delta_hu$ is bounded by $\Psi_+$.  After evaluating at $z'=z$, and
optimizing in
$\varepsilon$, this becomes
$$
\delta_hu(z,t')\leq \delta_hu(z,t)+ C(t'-t)^{{1+\alpha\over 2}}
+C|u|_{P^2(Q_1)}|t'-t|.
$$
A similar estimate from below follows by comparison
with a suitable subsolution.  The desired continuity in time of the first
spatial derivatives in $Q_{1/32}$ follows on sending $h$ to zero.

\subhead 3.4 Time regularity of second space derivatives\endsubhead

The proof of the parabolic $C^{2,\alpha}$ estimate can now be completed by
deducing appropriate continuity of the second spatial derivatives in time. 
We will deduce this as a consequence of the previous two estimates.

On $Q_{1/32}$ we have the estimates 
$$
|D^2u(x,t)-D^2u(y,t)|\leq C|x-y|^\alpha
$$
and
$$
|Du(x,t)-Du(x,t-\tau)|\leq C|\tau|^{{1+\alpha\over 2}}
$$
provided $t-\tau\geq -R^2/4$ and $x$, $y\in B_{R/2}$.
Fix $x\in B_{1/64}$ and $s,t\in(-1/256,0]$.  Let $\xi$ be an arbitrary
unit vector.  Then we have (provided $x+h\xi\in B_{1/32}$)
$$
D_\xi D_\xi u(x,t) \!= {D_\xi u(x+h\xi,t) - D_\xi u(x,t)\over h}
+{1\over h}\!\int_0^h\!\!(D_\xi D_\xi u(x+r\xi,t)-D_\xi D_\xi u(x,t))\,dr
$$
and so
$$
\eqalign{
|D_\xi D_\xi u(x,t)-D_\xi D_\xi u(x,s)|
&\leq h^{-1}\int_0^h\left|(D_\xi D_\xi u(x+r\xi,t)-D_\xi D_\xi
u(x,t)\right|\,dr\cr 
&\quad\null +h^{-1}\int_0^h\left|(D_\xi D_\xi
u(x+r\xi,s)-D_\xi D_\xi u(x,s)\right|\,dr\cr
&\quad\null+h^{-1}\left|D_\xi u(x+h\xi,t)-D_\xi u(x+h\xi,s)\right|\cr
&\quad\null +h^{-1}\left|D_\xi u(x,t)-D_\xi u(x,s)\right|\cr
&\quad\null\leq Ch^{\alpha}+Ch^{-1}|s-t|^{{1+\alpha\over 2}}.
\cr}
$$
Since $|s-t|< 1/256$ and $x\in B_{1/64}$, we can safely
choose $h=|s-t|^{1/2}< 1/64$ and still ensure that $x+h\xi\in B_{1/32}$.
This choice gives
$$
|D_\xi D_\xi u(x,t)-D_\xi D_\xi u(x,s)|\leq C|s-t|^{\alpha/2},
$$
proving the desired estimate.

\head 4.  Higher dimensions\endhead

In higher dimensions, a similar argument can be used to extend the
class of fully nonlinear parabolic equations for which second derivative
H\"older estimates hold:  Instead of concavity as a function of
the second derivatives, it suffices that the level sets of
$F(r,p,z,x,t)$ as a function of $r$ (for fixed 
$p$, $z$, $x$ and $t$) are convex. 

H\"older regularity of $u_t$ follows exactly as in Section 3.1. 
The key step is to show spatial regularity of the second space
derivatives in analogy with Section 3.2.  

The argument in Section 3.1 gives the estimate
$$
|u_t|_{P^{0,\alpha}(Q_{1/2})}\leq C\left(\text{\rm Lip}(F)
+\left|u_t\right|_{L^\infty(Q^1)}
\right).
$$
As in Section 3.2, for any fixed $t\in (-1/4,0]$ the uniformly elliptic
equation $G[D^2w(x),Dw(x),w(x),x]=\phi(x)$
holds, where $w(x)=u(x,t)$, $G[r,p,z,x]=F[r,p,z,x,t]$ and
$\phi(x)=u_t(x,t)$.  Note that
$\lambda\text{\rm I}\leq
\dot G\leq
\Lambda\text{\rm I}$ and
$\text{\rm Lip}(G)\leq\text{\rm Lip}(F)$.  Since the level sets of $F$ in
the first argument
are convex, and $F$ is uniformly monotone and $C^2$, there exists a constant
$K$ such that
$$
\ddot F^{klmn}M_{kl}M_{mn}\leq K\dot F^{kl}\dot F^{mn}M_{kl}M_{mn}
$$
for any symmetric matrix $M$.  But then $\tilde G=-\exp(-KG)$ is uniformly
monotone and concave in the first argument (with ellipticity constants
depending on $\sup_{Q_1}|D^2u|$), and we have $\tilde
G[D^2w,Dw,w,x]=\tilde\phi$ where $\tilde\phi=-\exp(-K\phi)$ is
H\"older-continuous on $B_{1/2}$. A perturbation result (see
\cite{C}, Theorem 3) then implies the estimate
$|D^2w(y)-D^2w(x)|\leq C|y-x|^{\alpha}$
for $x$ and $y$ in $B_{1/8}$.  The arguments of Sections 3.3
and 3.4 now apply unchanged to give the full result.

\enddocument

\end